\documentstyle[amscd,amssymb,verbatim,epsf]{amsart}
\begin{document}

\theoremstyle{plain}
\newtheorem{Thm}{Theorem}
\newtheorem{Cor}{Corollary}
\newtheorem{Con}{Conjecture}
\newtheorem{Main}{Main Theorem}
\newtheorem{Lem}{Lemma}
\newtheorem{Prop}{Proposition}
\newenvironment{Prf}{{\bf Proof:} }{\hfill $\Box$
\mbox{}}

\theoremstyle{definition}
\newtheorem{Def}{Definition}
\newtheorem{Note}{Note}

\theoremstyle{remark}
\newtheorem{notation}{Notation}
\renewcommand{\thenotation}{}

\errorcontextlines=0
\numberwithin{equation}{section}
\renewcommand{\rm}{\normalshape}%
\newcommand{\ici}[1]{\stackrel{\circ}{#1}}
\title[On Soft Mappings]%
   {On Soft Mappings}

\author{Cigdem Gunduz Aras* , Ayse Sonmez**, and H\"usey\.{I}n \c{C}akall\i\*** \\
Kocaeli University, Department of Mathematics, Kocaeli, Turkey Phone:(+90262)3032102\\ **Gebze Institute of Technology, Department of Mathematics,  Gebze-Kocaeli, Turkey Phone: (+90262)6051389 \\*** Maltepe University, Marmara E\u{g}\.{I}t\.{I}m K\"oy\"u, TR 34857, \.{I}stanbul-Turkey \; \; \; \; \; Phone:(+90216)6261050 ext:2248, \;  fax:(+90216)6261113 }

\address{Cigdem Gunduz Aras\\
                    Kocaeli University, Department of Mathematics, Kocaeli, Turkey Phone:(+90262)3032102}

\email{carasgunduz@@gmail.com; caras@@kocaeli.edu.tr}

\address{Ay\c{s}e S\"onmez\\
                     Department of Mathematics, Gebze Institute of Technology, Cayirova Campus 41400 Gebze- Kocaeli, Turkey Phone: (+90262)6051389}

\email{asonmez@@gyte.edu.tr; ayse.sonmz@@gmail.com}

\address{H\"usey\.{I}n \c{C}akall\i \\
          Maltepe University, Department of Mathematics, Marmara E\u{g}\.{I}t\.{I}m K\"oy\"u, TR 34857, Maltepe, \.{I}stanbul-Turkey \; \; \; \; \; Phone:(+90216)6261050 ext:2248, \;  fax:(+90216)6261113}

\email{hcakalli@@maltepe.edu.tr; hcakalli@@gmail.com}


\address{}

\email{}



\date{\today}

\maketitle

\begin{abstract}

In this paper, we introduce soft continuous mappings which are defined over an initial universe set with a fixed set of parameters. Later we study soft open and soft closed mappings, soft homeomorphism and investigate some properties of these concepts.

\end{abstract}

\maketitle

\section{Introduction}

\normalfont{}
Most of the real life problems in social sciences, engineering, medical sciences, economics etc.  the data involved are imprecise in nature. The solutions of such problems involve the use of mathematical principles based on uncertainty and imprecision. Thus classical set theory, which is based on the crisp and exact case may not be fully suitable for handling such problems of uncertainty.
A number of theories have been proposed for dealing with uncertainties efficiently way. Some of these are theory of fuzzy sets \cite{Zadeh[17]Fuzzysets}, theory of intuitionistic fuzzy sets \cite{Atanassov[3]Intuitionisticfuzzysets}, theory of vague sets, theory of interval mathematics \cite{Atanassov[4]Operatorsoverintervalvaluedintuitionisticfuzzysets} ,\cite{Gorzalzany[7]Amethodofinferenceinapproximatereasoningbasedoninterval-valuedfuzzysets} and theory of rough sets \cite{Pawlak[14]Roughsets}. However, these theories have their own difficulties. Molodtsov \cite{Molodtsov[10]Softsettheoryfirstresults} initiated a novel concept of soft sets theory as a new mathematical tool for dealing with uncertainties which is free from the above limitations. A soft set is a collection of approximate descriptions of an object. Soft systems provide a very general framework with the involvement of parameters. Since soft set theory has a rich potential, researchs on soft set theory and its applications in various fields are progressing rapidly.
Maji et al. \cite{MajiandBismasandRoy[11]Softsettheory} ,cite{MajiandRoyandBismas[12]AnApplicationofsoftsetsinadecisionmakingproblem}  worked on soft set theory and presented an application of soft sets in decision making problems.


M.Shabir and M.Naz \cite{ShabirandNaz[16]Onsofttopologicalspaces} introduced soft topological spaces. They also defined some concepts of soft sets on soft topological spaces. Later, researches about soft topological spaces were studied in \cite{ZorlutunaandAkdagandMinandAtmaca[18]Remarksonsofttopologicalspaces}-\cite{Min[20]Anoteonsofttopologicalspaces}.  In these studies, the concept of soft point is expressed by different approaches. In this paper we refer to the concept of soft point which was given in \cite{BayramovandAras[21]Softlocallycompactandsoftparacompactspaces} .\\

The purpose of this paper is to study soft continuity on soft topological spaces. Firstly, we recall some basic definitions about soft sets and the results from the literature. The continuity of mappings of soft topological spaces has defined and its properties has investigated. Finally soft open and soft closed mappings, soft homeomorphism are defined and some interesting results are derived which may be of value for further research.

\maketitle

\section{Preliminaries}

In this section, we give definitions and some results of soft sets.

\begin{Def} (\cite{Molodtsov[10]Softsettheoryfirstresults}) Let $X$ be an initial universe set and $E$ be a set of parameters. A pair $(F,E)$ is called a soft set over $X$ if only if $F$ is a mapping from $E$ into the set of all subsets of the set $X$, i.e., $F:E\rightarrow P(X)$, where $P(X)$ is the power set of $X$.
\end{Def}

\begin{Def}
(\cite{MajiandBismasandRoy[11]Softsettheory}) The intersection of two soft sets $(F,A)$ and $(G,B)$ over $X$ is the soft set $(H,C)$, where $C=A\cap B$ and $\forall c\in C$, $H(c)=F(c)\cap G(c)$. This is denoted by $(F,A)\tilde{\cap}(G,B)=(H,C)$.
\end{Def}

\begin{Def} (\cite{MajiandBismasandRoy[11]Softsettheory}) A soft set $(F,A)$ over $X$ is said to be a null soft set denoted by $\Phi$ if for all $a\in A$, $F(a)=\emptyset$ (null set).
\end{Def}

\begin{Def}(\cite{MajiandBismasandRoy[11]Softsettheory}) A soft set $(F,A)$ over $X$ is said to be an absolute soft set denoted by $\tilde{A}$ if for all $a\in A$, $F(a)=X$.
\end{Def}

\begin{Def}(\cite{MajiandBismasandRoy[11]Softsettheory}) The union of two soft sets $(F,A)$ and $(G,B)$ over $X$ is the soft set, where $C=A\cup B$ and $\forall c\in C$,\\
$H(\varepsilon)=\left\{
\begin{array}{c}
F(\varepsilon)\text{, if}\; \varepsilon \in A-B  \\
G(\varepsilon)\text{, if}\; \varepsilon \in B-A  \\
F(\varepsilon)\cup G(\varepsilon)\text{, if}\; \varepsilon \in A\cap B %
\end{array}%
\right. $.\\
This relationship is denoted by $(F,A)\tilde{\cup}(G,B)=(H,C)$.
\end{Def}

\begin{Def} (\cite{ShabirandNaz[16]Onsofttopologicalspaces} ) The complement of a soft set $(F,E)$ is denoted by $(F,E)^{\prime}$ and is defined by $(F,E)^{\prime}=(F^{\prime},E)$ where $F^{\prime}:E\rightarrow P(U)$ is a mapping given by $F^{\prime}(\alpha)=U-F(\alpha)$ for all $\alpha \in E$.
\end{Def}

\begin{Def}(\cite{ShabirandNaz[16]Onsofttopologicalspaces} ) Let $\tau$ be the collection of soft sets over $X$, then $\tau$ is said to be a soft topology on $X$ if\\
(1) $\Phi$, $\tilde{X}$  belong to $\tau$,\\
(2) the union of any number of soft sets in $\tau$ belongs to $\tau$,\\
(3) the intersection of any two soft sets in $\tau$ belongs to $\tau$.\\
The triplet $(X,\tau , E)$ is called a soft topological space over $X$.
\end{Def}

\begin{Def} (\cite{ShabirandNaz[16]Onsofttopologicalspaces} ) Let $(X,\tau , E)$ be a soft topological space over $X$. A soft set $(F,E)$ over $X$ is said to be a soft closed in $X$, if its relative complement $(F,E)^{\prime}$ belongs to $\tau$.
\end{Def}

\begin{Prop} (\cite{ShabirandNaz[16]Onsofttopologicalspaces} ) Let $(X,\tau , E)$ be a soft topological space over $X$. Then the collection $\tau_{\alpha}=\{F(\alpha)|(F,E)\in \tau\}$ for each $\alpha \in E$, defines a topology on $X$.
\end{Prop}

\begin{Def} (\cite{ShabirandNaz[16]Onsofttopologicalspaces} ) Let $(X,\tau , E)$ be a soft topological space over $X$ and $(F,E)$ be a soft set over $X$. Then the soft closure of $(F,E)$, denoted by $\overline{(F,E)}$ is the intersection of all soft closed super sets of $(F,E)$.
\end{Def}

\begin{Def} (\cite{ShabirandNaz[16]Onsofttopologicalspaces} ) Let $(X,\tau , E)$ be a soft topological space over $X$ and $(F,E)$ be a soft set over $X$. Then we associate with $(F,E)$ a soft set over $X$, denoted by $(\overline{F},E)$ and defined as $\overline{F}(\alpha)=\overline{F(\alpha)}$, where $\overline{F(\alpha)}$ is the closure of $F(\alpha)$ in $\tau_{\alpha}$ for each $\alpha \in E$.
\end{Def}

\begin{Prop} (\cite{ShabirandNaz[16]Onsofttopologicalspaces} ) Let $(X,\tau , E)$ be a soft topological space over $X$ and $(F,E)$ be a soft set over $X$. Then $(\overline{F},E)\subset \overline{(F,E)}$ .
\end{Prop}

\begin{Cor} (\cite{ShabirandNaz[16]Onsofttopologicalspaces}) Let $(X,\tau , E)$ be a soft topological space over $X$ and $(F,E)$ be a soft set over $X$. Then $(\overline{F},E)=\overline{(F,E)}$ if and only if $(\overline{F},E)^{\prime}\in \tau$.
\end{Cor}

\begin{Def} (\cite{HussainandAhmad[22]Somepropertiesofsofttopologicalspaces} ) Let $(X,\tau , E)$ be a soft topological space over $X$, $(G,E)$ be a soft set over $X$ and $x\in X$. Then $x$ is said to be a soft interior point of $(G,E)$, if there exists a soft open set $(F,E)$ such that $x\in (F,E)\subset (G,E)$.
\end{Def}

\begin{Def}(\cite{HussainandAhmad[22]Somepropertiesofsofttopologicalspaces}) Let $(X,\tau , E)$ be a soft topological space over $X$, $(G,E)$ be a soft set over $X$ and $x\in X$. Then $(G,E)$ is said to be a soft neighbourhood of $x$, if there exists a soft open set $(F,E)$ such that $x\in (F,E)\subset (G,E)$.
\end{Def}

\begin{Def}
(\cite{HussainandAhmad[22]Somepropertiesofsofttopologicalspaces}) Let $(X,\tau , E)$ be a soft topological space over $X$ then soft interior of soft set $(F,E)$ over $X$ is denoted by $ {(F,E)}^{\circ}$ and is defined as the union of all soft open sets contained in $(F,E)$.
\end{Def}
Thus $ {(F,E)}^{\circ}$ is the largest soft open set contained in $(F,E)$.

\maketitle

\section{Soft Continuous Mappings}

\begin{Def}(\cite{BayramovandAras[21]Softlocallycompactandsoftparacompactspaces} ) Let $(F,E)$ be a soft set over $X$. The soft set $(F,E)$ is called a soft point, denoted by $(x_e,E)$, if for the element $e\in E$, $F(e)=\{x\}$ and $F(e^{\prime})=\emptyset$ for all $e^{\prime}\in E-\{e\}$.
\end{Def}

\begin{Def} Let $(X,\tau , E)$ and $(Y,{\tau}^{\prime}, E)$ be two soft topological spaces, $f:(X,\tau , E)\rightarrow (Y,{\tau}^{\prime}, E)$ be a mapping. For each soft neighbourhood $(H,E)$ of $(f(x)_{e},E)$, if there exists a soft neighbourhood $(F,E)$ of $(x_{e},E)$ such that $f((F,E))\subset (H,E)$, then $f$ is said to be soft continuous mapping at $(x_{e},E)$.\\
If $f$ is soft continuous mapping for all $(x_{e},E)$, then $f$ is called soft continuous mapping.
\end{Def}

\begin{Thm} Let $(X,\tau , E)$ and $(Y,{\tau}^{\prime}, E)$ be two soft topological spaces, $f:(X,\tau , E)\rightarrow (Y,{\tau}^{\prime}, E)$ be a mapping. Then the following conditions are equivalent:\\
(1) $f:(X,\tau , E)\rightarrow (Y,{\tau}^{\prime}, E)$ is a soft continuous mapping,\\
(2) For each soft open set $(G,E)$ over $Y$, $f^{-1}((G,E))$ is a soft open set over $X$,\\
(3) For each soft closed set $(H,E)$ over $Y$, $f^{-1}((H,E))$ is a soft closed set over $X$,\\
(4) For each soft set $(F,E)$ over $X$, $f(\overline{(F,E)})\subset \overline{(f(F,E))}$,\\
(5) For each soft set $(G,E)$ over $Y$, $\overline{(f^{-1}(G,E))}\subset f^{-1}(\overline{(G,E)})$,\\
(6) For each soft set $(G,E)$ over $Y$, $f^{-1}((G,E)^{\circ})\subset (f^{-1}(G,E))^{\circ}$.\\
\end{Thm}
\begin{Prf} $(1)\Rightarrow (2)$ Let $(G,E)$ be a soft open set over $Y$ and $(x_{e},E)\in f^{-1}(G,E)$ be an arbitrary soft point. Then $f(x_{e},E)=(f(x)_{e},E)\in (G,E)$. Since $f$ is soft continuous mapping, there exists $(x_{e},E)\in (F,E)\in \tau$ such that $f(F,E)\subset (G,E)$. This implies that $(x_{e},E)\in (F,E)\subset f^{-1}(G,E)$, $f^{-1}((G,E))$ is a soft open set over $X$.\\

$(2)\Rightarrow (1)$ Let $(x_{e},E)$ be a soft point and $(f(x)_{e},E)\in (G,E)$ be an arbitrary soft neighbourhood. Then $(x_{e},E)\in f^{-1}(G,E)$ is a soft neighbourhood and $f(f^{-1}(G,E))\subset (G,E)$.\\

$(3)\Rightarrow (4)$ Let $(F,E)$ be a soft set over $X$. Since $(F,E)\subset f^{-1}(f(F,E))$ and $f(F,E)\subset \overline{(f(F,E))}$, we have $(F,E)\subset f^{-1}(f(F,E))\subset f^{-1}\overline{(f(F,E))}$. By part (3), since $f^{-1}\overline{(f(F,E))}$ is a soft closed set over $X$, $\overline{(F,E)}\subset f^{-1}\overline{(f(F,E))}$. Thus $f(\overline{(F,E)})\subset f(f^{-1}\overline{(f(F,E))})\subset \overline{f(F,E)}$ is obtained.\\

$(4)\Rightarrow (5)$ Let $(G,E)$ be a soft set over $Y$ and $f^{-1}(G,E)=(F,E)$. By part (4), we have $f(\overline{(F,E)})=f(\overline{f^{-1}(G,E)})\subset \overline{f(f^{-1}(G,E))}\subset \overline{(G,E)}$. Then $\overline{f^{-1}(G,E)}=\overline{(F,E)}\subset f^{-1}(f\overline{(F,E)})\subset f^{-1}(\overline{(G,E)})$.\\

$(5)\Rightarrow (6)$ Let $(G,E)$ be a soft set over $Y$. Substituting $(G,E)^{\prime}$ for condition in (5). Then $\overline{f^{-1}((G,E)^{\prime})}\subset f^{-1}(\overline{(G,E)^{\prime}})$. Since $(G,E)^{\circ}=(\overline{(G,E)^{\prime}})^{\prime}$, then we have $f^{-1}((G,E)^{\circ})=f^{-1}((\overline{(G,E)^{\prime}})^{\prime})=(f^{-1}(\overline{(G,E)^{\prime}}))^{\prime}\subset (\overline{f^{-1}((G,E)^{\prime})})^{\prime}=(\overline{(f^{-1}(G,E))^{\prime}})^{\prime}=(f^{-1}(G,E))^{\circ}$.\\

$(6)\Rightarrow (2)$ Let $(G,E)$ be a soft open set over $Y$. Then since $(f^{-1}(G,E))^{\circ}\subset f^{-1}(G,E)=f^{-1}((G,E)^{\circ})\subset (f^{-1}(G,E))^{\circ}$, $(f^{-1}(G,E))^{\circ}=f^{-1}(G,E)$ is obtained. This implies that $f^{-1}(G,E)$ is a soft open set over $X$.
\end{Prf}

\textbf{Example 1.} Let $X=\{h_1,h_2,h_3\}$, $E=\{e_1,e_2\}$ and $\tau =\{\Phi , \tilde{X}, (F_1,E), (F_2,E)\}$, ${\tau}^{\prime}=\{\Phi , \tilde{X}, (G_1,E), (G_2,E)\}$ be two soft topologies defined on $X$, where $(F_1,E)$, $(F_2,E)$, $(G_1,E)$ and $(G_2,E)$ are soft sets over $X$, defined as follows:
\[F_1(e_1)=\{h_1,h_2\}, \; \; \;F_1(e_2)=\{h_3\}, \; \; \;F_2(e_1)=X, \; \; \; F_2(e_2)=\{h_3\}, \; \; \;\]
and
\[G_1(e_1)=\{h_1\}, \; \; \;G_1(e_2)=\{h_3\}, \; \; \;G_2(e_1)=\{h_1,h_3\}, \; \; \;G_2(e_2)=\{h_2,h_3\}, \; \; \;\]
If we get the mapping $f:X\rightarrow X$ defined as
\[f(h_1)=f(h_2)=h_1, \; f(h_3)=h_3\]
then since $f^{-1}(G_1,E)=(F_1,E)$ and $f^{-1}(G_2,E)=(F_2,E)$, $f$ is a soft continuous mapping.\\

\textbf{Example 2.} Let $X=\{h_1,h_2,h_3\}$, $E=\{e_1,e_2\}$ and $\tau=\{\Phi ,\tilde{X}, (F_1,E), (F_2,E), (F_3,E), (F_4,E) \}$, ${\tau}^{\prime}=\{\Phi ,\tilde{X}, (G_1,E), (G_2,E), (G_3,E), (G_4,E) \}$ be two soft topologies defined on $X$ where $(F_1,E), (F_2,E), (F_3,E), (F_4,E), (G_1,E), (G_2,E), (G_3,E)$ and $(G_4,E)$ are soft sets over $X$, defined as follows:
\[F_1(e_1)=\{h_2\}, \; \; \;F_1(e_2)=\{h_1\}, \; \; \;F_2(e_1)=\{h_2,h_3\}, \; \; \; F_2(e_2)=\{h_1,h_2\}, \; \; \;\]
\[F_3(e_1)=\{h_3\}, \; \; \;F_3(e_2)=\{h_1,h_2\}, \; \; \;F_4(e_1)=\emptyset, \; \; \; F_4(e_2)=\{h_1\}, \; \; \;\]
\[F_5(e_1)=X, \; \; \; F_2(e_2)=\{h_1,h_2\}\]
and
\[G_1(e_1)=\{h_2\}, \; \; \;G_1(e_2)=\{h_1\}, \; \; \;G_2(e_1)=\{h_2,h_3\}, \; \; \;G_2(e_2)=\{h_1,h_2\}, \; \; \;\]
\[G_3(e_1)=\{h_1,h_2\}, \; \; \;G_3(e_2)=X, \; \; \;G_4(e_1)=\{h_2\}, \; \; \;G_4(e_2)=\{h_1,h_2\}, \; \; \;\]
Then $(X, \tau ,E)$, $(X, {\tau}^{\prime}, E)$ are two soft topological spaces and $f=1_X:X\rightarrow X$ is not soft continuous mapping.

\begin{Thm} If $f:(X, \tau ,E)\rightarrow (Y, {\tau}^{\prime}, E)$ is a soft continuous mapping, then for each $\alpha \in E$, $f_{\alpha}:(Y,{\tau}_{\alpha})\rightarrow (Y,{\tau}_{\alpha}^{\prime})$ is a continuous mapping.
\end{Thm}
\begin{Prf} Let $U\in {\tau}_{\alpha}^{\prime}$. Then there exists a soft open set  $(G,E)$ over $Y$ such that $U=G(\alpha)$. Since $f:(X, \tau ,E)\rightarrow (Y, {\tau}^{\prime}, E)$ is a soft continuous mapping, $f^{-1}(G,E)$ is a soft open set over $X$ and $f^{-1}(G,E)(\alpha)=f^{-1}(G(\alpha))=f^{-1}(U)$ is an open set. This implies that $f_{\alpha}$ is a continuous mapping.
\end{Prf}

Now we give an example to show that the converse of above theorem does not hold.

\textbf{Example 3.} Let $X=\{x_1,x_2,x_3\}$, $Y=\{y_1,y_2,y_3\}$ and $E=\{e_1,e_2\}$. Then $\tau=\{\Phi ,\tilde{X}, (F_1,E), (F_2,E), (F_3,E), (F_4,E), (F_5,E) \}$ is a soft topological space over $X$ and ${\tau}^{\prime}=\{\Phi ,\tilde{Y}, (G_1,E), (G_2,E), (G_3,E) \}$ is a soft topological space over $Y$. Here $(F_1,E), (F_2,E), (F_3,E), (F_4,E), (F_5,E)$ are soft sets over $X$ and   $(G_1,E), (G_2,E), (G_3,E)$ are soft sets over $Y$, defined as follows:

\[F_1(e_1)=\{x_1\}, \; \; \;F_1(e_2)=\{x_1,x_3\}, \; \; \; F_2(e_1)=\{x_2\}, \; \; F_2(e_2)=\{x_1\} ,\]
\[F_3(e_1)=\{x_1,x_2\}, \; \; F_3(e_2)=\{x_1,x_3\},\; \; \; F_4(e_1)=\emptyset, \; \; \;F_4(e_2)=\{x_1\},  \; \; \; \;F_5(e_1)=\{x_1,x_2\}, \; \; \; F_5(e_2)=X.\]
and
\[G_1(e_1)=Y, \; \; \; G_1(e_2)=\{y_2\}, \; \; \;G_2(e_1)=\{y_1\}, \; \; \; G_2(e_2)=\{y_2\},\; \; \; G_3(e_1)=\{y_1,y_2\}, \; \; \; G_3(e_2)=\{y_2\}.\]

If we get the mapping $f:X\rightarrow Y$ defined as
\[f(x_1)=y_2, \; f(x_2)=y_1, \; f(x_3)=y_3\]
then $f$ is not a soft continuous mapping, since $f^{-1}(G_1)\notin \tau$, where $f^{-1}(G_1)(e_1)=X$, $f^{-1}(G_1)(e_2)=\{x_1\}$. Also, $f_{e_1}:(X, \tau_{e_1})\rightarrow (Y,\tau_{e_1}^{\prime})$ and $f_{e_2}:(X, \tau_{e_2})\rightarrow (Y,\tau_{e_2}^{\prime})$ are continuous mappings. Here
\[\tau_{e_1}=\{\emptyset , X, \{x_1\}, \{x_2\}, \{x_1,x_2\}\}, \; \; \; \tau_{e_2}=\{\emptyset , X, \{x_1\},  \{x_1,x_3\}\}\]
and
\[\tau_{e_1}^{\prime}=\{\emptyset , Y, \{y_1\}, \{y_1,y_2\}\}, \; \; \; \tau_{e_2}^{\prime}=\{\emptyset , Y, \{y_2\}\}.\]

Now, let us show that when the above theorem is true.

\begin{Thm} If $(\overline{F},E)^{\prime}$ is a soft open set over $X$, for each soft set $(F,E)$, then $f:(X, \tau ,E)\rightarrow (Y, {\tau}^{\prime}, E)$  is a soft continuous mapping if and only if $f_{\alpha}:(X, \tau_{\alpha})\rightarrow (Y, {\tau}^{\prime}_{\alpha})$  is continuous mapping, for each $\alpha \in E$.
\end{Thm}
\begin{Prf} Let $f_{\alpha}:(X, \tau_{\alpha})\rightarrow (Y, {\tau}^{\prime}_{\alpha})$ be a continuous mapping, for each $\alpha \in E$, and let $(F,E)$ be an arbitrary soft set over $X$. Then $f_{\alpha}(\overline{F}(\alpha))\subset \overline{f(F)}(\alpha)$ is satisfied, for each $\alpha \in E$. Since $(\overline{F},E)^{\prime}\in \tau$, $(\overline{F},E)=\overline{(F,E)}$ from Corollary 1. Thus
$f(\overline{(F,E)})\subset \overline{f((F,E))}$ is obtained. This implies that $f:(X, \tau ,E)\rightarrow (Y, {\tau}^{\prime}, E)$ is a soft continuous mapping.
\end{Prf}

\begin{Def} Let $(X, \tau ,E)$ and $(Y, {\tau}^{\prime}, E)$ be two soft topological spaces, $f:X\rightarrow Y$ be a mapping.\\
a) If the image $f((F,E))$ of each soft open set $(F,E)$ over $X$ is a soft open set in $Y$, then $f$ is said to be a soft open mapping.\\
b) If the image $f((H,E))$ of each soft closed set $(H,E)$ over $X$ is a soft closed set in $Y$, then $f$ is said to be a soft closed mapping.
\end{Def}

\begin{Prop} If $f:(X, \tau ,E)\rightarrow (Y, {\tau}^{\prime}, E)$ is soft open (closed), then for each $\alpha \in E$, $f_{\alpha}:(X, \tau_{\alpha})\rightarrow (Y, {\tau}^{\prime}_{\alpha})$ is an open (closed) mapping.
\end{Prop}

\begin{Prf}The proof of the proposition is straightforward and it is left to the reader.
\end{Prf}

Note that the concepts of soft continuous, soft open and soft closed mappings are all independent of each other.

\textbf{Example 4.} Let $(X,\tau ,E)$ be soft discrete topological space and $(X,\tau^{\prime} ,E)$ be soft indiscrete topological space. Then $1_X:(X,\tau ,E)\rightarrow (X,\tau^{\prime} ,E)$ is a soft open and soft closed mapping. But it is not soft continuous mapping.\\

\textbf{Example 5.} Let $X=\{h_1,h_2,h_3\}$, $E=\{e_1,e_2\}$ and $\tau=\{\Phi ,\tilde{X}, (F_1,E), (F_2,E),..., (F_7,E) \}$, ${\tau}^{\prime}=\{\Phi ,\tilde{X}, (G_1,E), (G_2,E), (G_3,E), (G_4,E) \}$ be two soft topologies defined on $X$ where $(F_1,E), (F_2,E), (F_3,E),..., (F_7,E), (G_1,E), (G_2,E), (G_3,E)$  and $(G_4,E)$ are soft sets over $X$, defined as follows:
\[F_1(e_1)=\{h_2\}, \; \; \;F_1(e_2)=\{h_1\}, \; \;  \; \;F_2(e_1)=\{h_1,h_3\}, \; \; \; F_2(e_2)=\{h_2,h_3\},\]
\[F_3(e_1)=\{h_2\}, \; \; \;F_3(e_2)=X, \; \; \;  \;F_4(e_1)=\emptyset , \; \; \; F_4(e_2)=\{h_1\},\]
\[F_5(e_1)=\{h_1,h_3\}, \; \; \;F_5(e_2)=X, \; \; \;  \;F_6(e_1)=\emptyset , \; \; \; F_6(e_2)=\{h_2,h_3\},\]
\[F_7(e_1)=\emptyset,  \; \; \; F_7(e_2)=X \]
and
\[G_1(e_1)=\{h_2\}, \; \; \;G_1(e_2)=\{h_1\}, \; \;  \; \;G_2(e_1)=\{h_2,h_3\}, \; \; \; G_2(e_2)=\{h_1,h_2\},\]
\[G_3(e_1)=\{h_1,h_2\}, \; \; \;G_3(e_2)=X, \; \;  \; \;G_4(e_1)=\{h_2\}, \; \; \; G_4(e_2)=\{h_1,h_2\}.\]
If we get the mapping $f:X\rightarrow X$ defined as $f(h_i)=h_1$, for $1\leq i\leq 3$. It is clear that
\[f^{-1}(G_1)(e_1)=f^{-1}(G_4)(e_1)=\emptyset , \; \; f^{-1}(G_1)(e_2)=f^{-1}(G_4)(e_2)=X, \; \; \; f^{-1}(G_3)(e_1)=X, \; \; f^{-1}(G_3)(e_2)=X.\]
Then $f$ is a soft continuous mapping, but
\[f(F_1)(e_1)=\{h_1\}, \; \; f(F_1)(e_2)=\{h_1\}, \; \; f(F_1^{\prime})(e_1)=\{h_1\}, \; \; f(F_1^{\prime})(e_2)=\{h_1\}.\]
Hence it is not both soft open and soft closed mapping.

\textbf{Example 6.} Let $X=\{h_1,h_2,h_3\}$, $Y=\{a,b\}$ and $E=\{e_1,e_2\}$ and $\tau=\{\Phi ,\tilde{X}, (F_1,E), (F_2,E)\}$, ${\tau}^{\prime}=\{\Phi ,\tilde{Y}, (G_1,E), (G_2,E)\}$ be two soft topologies defined on $X$ and $Y$, respectively. Here $(F_1,E), (F_2,E), (G_1,E), (G_2,E)$ are soft sets over $X$  and $Y$, respectively. The soft sets are defined as follows:
\[F_1(e_1)=\{h_1,h_2\}, \; \; \;F_1(e_2)=\{h_3\}, \; \; \; \;F_2(e_1)=X, \; \; \; F_2(e_2)=\{h_3\},\]
and
\[G_1(e_1)=Y, \; \; \;G_1(e_2)=\{b\},  \; \; \;G_2(e_1)=\{a\}, \; \; \; G_2(e_2)=\{b\},\]
If we get the mapping $f:X\rightarrow Y$ defined as
\[f(h_1)=\{a\},\; \; \;  f(h_2)=f(h_3)=\{b\}.\]
It is clear that
\[f(F_1)(e_1)=Y, \; \; f(F_1)(e_2)=\{b\}, \; \; f(F_2)(e_1)=Y, \; \; f(F_2)(e_2)=\{b\}.\]
Then the mapping $f:X\rightarrow Y$ is a soft open mapping. Also since $f(F_1^{\prime})(e_1)=\{b\}, \; \; f(F_1^{\prime})(e_2)=Y$, it is not soft closed mapping and $f^{-1}(G_1)(e_1)=X$, $f^{-1}(G_1)(e_2)=\{h_2,h_3\}$. Hence it is not soft continuous mapping.

\textbf{Example 7.} Let  $X=\{h_1,h_2,h_3\}$, $Y=\{a,b\}$, $E=\{e_1,e_2\}$ and $\tau=\{\Phi ,\tilde{X}, (F_1,E), (F_2,E), (F_3,E)\}$, ${\tau}^{\prime}=\{\Phi ,\tilde{Y}, (G_1,E), (G_2,E)\}$ be two soft topologies defined on $X$ and $Y$, respectively. Here $(F_1,E), (F_2,E), (F_3,E), (G_1,E), (G_2,E)$ are soft sets over $X$ and $Y$, respectively. The soft sets are defined as follows:
\[F_1(e_1)=\{h_1,h_3\},  \;F_1(e_2)=\{h_2\},  \;F_2(e_1)=X,  \; F_2(e_2)=\{h_2,h_3\},  \; F_3(e_1)=\{h_3\},  \; F_3(e_2)=\{h_2\}\]
and
\[G_1(e_1)=\emptyset, \; \; \;G_1(e_2)=\{a\},  \; \; \;G_2(e_1)=\{a\}, \; \; \; G_2(e_2)=Y.\]
Now we define the mapping $f:X\rightarrow Y$ as $f(h_1)=f(h_2)=\{a\}$, $f(h_3)=\{b\}$. It is clear that $f(F_1^{\prime}(e_1))=f(h_2)=\{a\}, f(F_1^{\prime}(e_2))=f(\{h_1,h_3\})=Y$, $f(F_2^{\prime}(e_1))=\emptyset, f(F_2^{\prime}(e_2))=f(\{h_1\})=\{a\}$, $f(F_3^{\prime}(e_1))=\{a\}, f(F_3^{\prime}(e_2))=Y$. This implies that $f$ is a soft closed mapping.
Also $f(F_1(e_1))=Y, f(F_1(e_2))=\{a\}$, $f^{-1}(G_1(e_1))=\emptyset, f^{-1}(G_1(e_2))=\{h_1,h_2\}$. Then it is not soft open and soft continuous mapping, respectively.

\begin{Thm} Let $(X,\tau, E)$ and $(Y,\tau^{\prime},E)$ be two soft topological spaces, $f:X\rightarrow Y$ be a mapping.\\
a) $f$ is a soft open mapping if and only if for each soft set $(F,E)$ over $X$, $f((F,E)^{\circ})\subset (f(F,E))^{\circ}$ is satisfied.\\
b) $f$ is a soft closed mapping if and only if for each soft set $(F,E)$ over $X$, $(\overline{f(F,E))}\subset f(\overline{(F,E)})$ is satisfied.
\end{Thm}
\begin{Prf} a) Let $f$ be a soft open mapping and $(F,E)$ be a soft set over $X$. $(F,E)^{\circ}$ is a soft open set and $(F,E)^{\circ}\subset (F,E)$. Since $f$ is a soft open mapping, $f((F,E)^{\circ})$ is a soft open set in $Y$ and $f((F,E)^{\circ})\subset f((F,E))$ . Thus $f((F,E)^{\circ})\subset f((F,E))^{\circ}$ is obtained.\\
Conversely, let $(F,E)$ be any soft open set over $X$. Then $(F,E)=(F,E)^{\circ}$. From the condition of theorem, we have $f((F,E)^{\circ})\subset (f(F,E))^{\circ}$. Then $f((F,E))=f((F,E)^{\circ})\subset (f(F,E))^{\circ}\subset f((F,E))$. This implies that $f((F,E))=(f(F,E))^{\circ}$. This completes the proof.\\
b) Let $f$ be a soft closed mapping and $(F,E)$ be any soft set over $X$. Since $f$ is a soft closed mapping, $f(\overline{(F,E)})$ is a soft closed set over $Y$ and $f((F,E))\subset f(\overline{(F,E)})$. Thus $\overline{f(F,E)}\subset f(\overline{(F,E)})$ is obtained.\\
Conversely, let $(F,E)$ be any soft closed set over $X$. From the condition of theorem, $\overline{(f(F,E))}\subset f(\overline{(F,E)})=f((F,E))\subset \overline{(f(F,E))}$. This means that $\overline{(f(F,E))}=f((F,E))$. This completes the proof.
\end{Prf}
\begin{Def} Let $(X,\tau, E)$ and $(Y,\tau^{\prime},E)$ be two soft topological spaces, $f:X\rightarrow Y$ be a mapping. If $f$ is a bijection, soft continuous and $f^{-1}$ is a soft continuous mapping, then $f$ is said to be soft homeomorphism from $X$ to $Y$. When a homeomorphism $f$ exists between $X$ and $Y$, we say that $X$ is soft homeomorphic to $Y$.
\end{Def}
\begin{Thm} Let $(X,\tau, E)$ and $(Y,\tau^{\prime},E)$ be two soft topological spaces, $f:X\rightarrow Y$ be a bijective mapping. Then the following conditions are equivalent:\\
(1) $f$ is a soft homeomorphism,\\
(2) $f$ is a soft continuous and soft closed mapping,\\
(3) $f$ is a soft continuous and soft open mapping.
\end{Thm}
\begin{Prf} It is easily obtained.
\end{Prf}

\section{Conclusion}
We have introduced soft continuous mappings which are defined over an initial universe with a fixed set of parameters. Later we study soft open and soft closed mappings, soft homeomorphism and investigate some properties of these concepts. In the end, we must say that, soft topological spaces are defined over different set of parameters and it does not make any sense for the results of this paper.

\end{document}